# La leçon mathématique de Théodore (*Théétète*, 147c7-d6) I

## Mais que dessine donc Théodore ?

Salomon Ofman

### I.  Introduction

Les deux exposés portent sur la première partie de ce qu'on nomme le 'passage mathématique du *Théétète*' (147c7-d6), où Théétète raconte brièvement la leçon de mathématique donnée par le géomètre Théodore. Le premier porte sur les dessins tracés par Théodore, le second sur ses méthodes de démonstration d'incommensurabilité.

Ces exposés sont basés sur un travail en commun avec Luc Brisson, commencé il y a de nombreuses années. Nous l'avons repris il y a quelques mois pour un livre en cours de rédaction sur ce passage. Celui-ci a été commenté dès l'Antiquité et a donné lieu à un très grand nombre de commentaires des philosophes, des historiens et même des mathématiciens. Les chercheurs modernes ne sont d'accord à peu près sur rien, si ce n'est que son sujet est l'incommensurabilité de certaines grandeurs géométriques.

La meilleure interprétation de ce passage est probablement celle de Wilbur Knorr dans son ouvrage *The Evolution of Euclidean* Elements (ou *EEE*). Mais elle pose de nombreuses difficultés. Ce qui est important est que Knorr ne se contente pas comme la plupart des interprètes de montrer en détails sa reconstruction, mais il donne six critères que doit vérifier toute interprétation pour être crédible (cf. Annexe II). Comme il le remarque lui-même, cela ouvre la voie à un nombre indéfini d'autres méthodes qui seraient alors aussi voire plus crédibles que la sienne. Heureusement ajoute-t-il, et à raison, 'nous verrons qu'en fait, aucune reconstruction (avant celle présentée ici [dans *EEE*]) ne les vérifie.' (*EEE*, p. 96-7)

Nous sommes d'accord avec tout ce qui est dit ici. Mais alors, pourquoi présenter une nouvelle interprétation ?

C'est que, si selon nous ces six conditions sont nécessaires, il faut en ajouter trois autres essentielles, dont aucune n'est vérifiée par les différentes reconstructions, celle de Knorr y compris :

1) On ne doit pas postuler *a priori* que Platon veut célébrer ici Théétète, futur grand mathématicien, et qu'il lui rend hommage ici.

2) Le passage doit être lu comme partie du *Théétète* et non pas comme une partie autonome.

3) Sa description par Platon est réaliste au sens où elle doit suivre de près les leçons réelles des mathématiciens sur le sujet, que cette leçon soit mise en scène par Platon ou ait eu vraiment lieu.

Celle que nous allons donner vérifie alors non seulement les six contraintes posées par Knorr, mais les trois autres ci-dessus.

### II.  Le texte

Nous faisons quelques changements par rapport à celui communément utilisé. Alors que les interprétations débutent généralement en 147d3, nous commençons 3 lignes plus haut comme vous l'avez vu dans le titre. D'autre part, nous apportons deux modifications à l'édition classique de Burnet dont une ponctuation : comme vous le savez, les textes de l'époque de Platon ne comportaient ni ponctuation, ni minuscules, ni accents.

| Théétète | Θεαίτητος |
|---|---|
| Et en fait, il se pourrait que tu demandes la sorte de chose qui nous est venue à l'esprit, à moi et à ce Socrate-là ton homonyme, alors que nous discutions entre nous tout à l'heure. | (147c7) ἀτὰρ κινδυνεύεις ἐρωτᾶν οἷον καὶ αὐτοῖς ἡμῖν ἔναγχος εἰσῆλθε διαλεγομένοις, ἐμοί τε καὶ τῷ σῷ ὁμωνύμῳ τούτῳ Σωκράτει. |
| **Socrate** | **Σωκράτης** |
| Quoi donc, Théétète ? | τὸ ποῖον δή, ὦ Θεαίτητε; |
| **Théétète** | **Θεαίτητος** |
| Quelque chose à propos des puissances. Théodore ici présent nous avait trace des figures, mettant en évidence que la puissance de 3 pieds aussi bien que celle de 5 pieds ne sont pas commensurables en longueur avec celle de 1 pied ; il a continué ainsi, prenant chacune à tour de rôle, jusqu'à celle de 17 pieds. Il s'est arrêté là pour quelque raison. | (147d2) περὶ δυνάμεών τι. ἡμῖν Θεόδωρος ὅδε ἔγραφε, τῆς τε τρίποδος πέρι καὶ πεντέποδος ἀποφαίνων ὅτι μήκει οὐ σύμμετροι τῇ ποδιαίᾳ, καὶ οὕτω κατὰ μίαν ἑκάστην προαιρούμενος μέχρι τῆς ἑπτακαιδεκάποδος· ἐν δὲ ταύτῃ πως ἐνέσχετο. |

Modifications par rapport à l'édition de Burnet :

- Un point après le 'τι' de la première ligne dans la deuxième intervention de Théétète.
- L'acceptation de 'ἀποφαίνων' qui est mis entre parenthèses par Burnet.

**Explication** : Socrate vient juste de reprocher à Théétète de ne pas comprendre ce qu'est une définition, car pour définir la 'science', il a donné une liste de 'τέχναι' (146d1). Socrate lui donne alors l'exemple de la définition de la glaise ('πηλός') mélange d'eau et de terre. Sur quoi, Théétète déclare que maintenant c'est facile ('ῥᾴδιον'), et dévie sur une discussion qu'il a eue avec un camarade homonyme du philosophe. Socrate lui demande l'objet de cette discussion, et Théétète lui répond qu'il s'agit de 'quelque chose' à propos des '*dynameôn*'. Il explique alors que Théodore leur avait dessiné des figures dans une leçon concernant des incommensurables.

Dans les lectures usuelles, on comprend que le jeune garçon affirme directement que 'Théodore a dessiné quelque chose' à leur propos, ce qui n'est pas une réponse à la question et n'ajoute rien au reste du texte. En fait, cela n'a guère de sens de dire que l'on dessine 'à propos' de quelque chose. On dessine non pas 'à propos' des triangles, des cercles ou des carrés, mais des cercles, des triangles ou des carrés. Si l'on veut donner un sens à la phrase dans le cadre de l'interprétation usuelle, on ne peut pas comprendre 'ἔγραφε' comme signifiant simplement 'dessiner', mais plutôt par 'prouver par des dessins (ou des diagrammes)'.

Au contraire, selon notre interprétation, il est deux temps. D'une part, Théodore *dessine* les 'puissances', puis dans un second temps, il prouve qu'elles sont incommensurables 'en longueur' à l'unité. Cela explique sans doute pourquoi Burnet décide de supprimer 'ἀποφαίνων' qui fait effectivement double emploi dans l'interprétation usuelles, dont la sienne.

Mais alors qu'a dessiné Théodore ?

### III.    Les dessins de Théodore

L'interprétation de '*dynamis/eis*' a été l'objet d'innombrables discussions pour savoir s'il s'agit de carrés ou des côtés du carré. Je n'aborderai pas cette question difficile qui ne joue pas de rôle important ici. Commençons tout d'abord par un rappel sur la définition de grandeurs incommensurables qui jouent un rôle essentiel dans le passage mathématique.

**Commensurabilité/incommensurabilité** : deux grandeurs géométriques, par exemple deux longueurs *A* et *B*, sont dites commensurables entre elles s'il existe une 'commune mesure', c'est-à-dire une unité *U* telle que à la fois *A* et *B* soient des multiples entiers de *U*. En écriture moderne : $A = aU$ et $B = bU$ où *a* et *b* sont des entiers. Si une telle mesure *U* n'existe pas, les grandeurs *A* et *B* sont incommensurables.

Le lien avec la rationalité et l'irrationalité est donné en remarquant que *A* et *B* sont commensurables si et seulement si leur rapport $A/B = a/b$ est égal au rapport de 2 entiers *a* et *b*, c'est-à-dire si et seulement si leur rapport $A/B$ est rationnel. Inversement, *A* et *B* sont incommensurables si et seulement si leur rapport $A/B$ est irrationnel.

**Exemples** : Tout d'abord, deux grandeurs entières sont toujours commensurables, puisque toutes deux sont des multiples entiers de l'unité. Ainsi les longueurs de *3* pieds et *5* pieds sont commensurables. Et leur rapport *3/5* est rationnel. Mais il n'est pas nécessaire que les grandeurs soient entières pour être commensurables entre elles.

Ainsi soient *A* et *B* des longueurs respectivement de *3* pieds et *5/2* pieds. Elles sont commensurables car en prenant pour unité *U = 1/2* pied, on a : $A = 6U$ et $B = 5U$. Et le rapport $A/B = 6/5$ est bien rationnel.

Par contre, d'après le théorème de Pythagore, la longueur de la diagonale D du carré de côté *2* pieds est (en écriture moderne) égale à $2\sqrt{2}$ pieds. En effet :

$D^2 = 2^2 + 2^2 = 4 + 4 = 8$ pieds, d'où $D = 2\sqrt{2}$ pieds.

Pour les géomètres grecs, cela revient à dire que le carré construit sur la diagonale du carré de *2* pieds est le double de ce carré, c'est ce que montre Socrate au jeune esclave dans le *Ménon*.

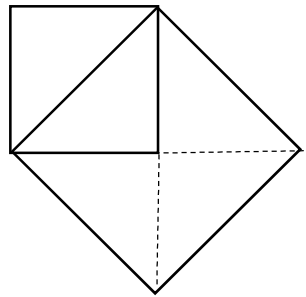

**Figure 1**

Le rapport de la diagonale au côté du carré est donc : $2\sqrt{2}/2 = \sqrt{2}$. Leur rapport est irrationnel, elle est donc incommensurable au côté du carré.

**Le problème du pied ('πούς') comme unité**

Le 'pied' ('πούς') dont il est question ici est d'abord une unité usuelle de longueur en Grèce ancienne. Mais elle est aussi utilisée ici par Théétète comme unité de surface, comme le fait Socrate dans le *Ménon* lorsqu'il montre le doublement du carré au jeune esclave. En termes modernes, on parlerait d'un pied-carré. Les surfaces étant des nombres entiers de pieds, elles sont toujours commensurables avec celle de 1 pied. Théétète explique à Socrate que Théodore leur a dessiné des grandeurs qui sont incommensurables avec celle de longueur 1 pied. Il s'agit donc des côtés de carrés de surfaces entières et le pied auquel sont comparés ces côtés est cette fois la longueur de 1 pied.

**La suite de Théodore**

Théodore dessine donc tout d'abord les côtés des carrés de *3* pieds, puis de *5* pieds, et il continue ainsi jusqu'à celui de *17* pieds. Il y a de nombreuses discussions sur leur construction, une question très débattue est de savoir si le carré de *17* pieds est compris ou pas dans la suite étudiée par Théodore. Je considérerai cette question dans le deuxième exposé.

Quoiqu'il en soit Théétète dit seulement que cette suite est la suivante : *3, 5, …, 17* (compris ou pas), sans préciser plus. La thèse quasi-universellement acceptée est qu'il s'agit de la suite des entiers compris entre *3* et *17* (*17* compris ou pas) excepté les carrés parfaits, soit :

**Suite standard de Théodore** : ***3, 5, 6, 7, 8, 10, 11, 12, 13, 14, 15, (17)***.

L'absence de *2* est alors justifiée par la référence à l'incommensurabilité de la diagonale par rapport au côté du carré. En effet, comme Socrate le montre dans le *Ménon*, la surface du carré de côté la diagonale d'un carré donné est le double du carré original. Partant du carré de 1 pied, la diagonale construite sur ce carré correspond alors au cas du carré de *2* pieds, et son absence s'explique par l'ancienneté de la connaissance de leur incommensurabilité. Cela est certainement juste. Toutefois, la preuve de cette incommensurabilité n'est jamais prise en compte dans les interprétations usuelles, alors que nous allons voir qu'elle a un rôle essentiel pour rendre compte des grandeurs considérées par Théodore, et ensuite dans les démonstrations qui sont l'objet du second exposé.

**L'incommensurabilité de la diagonale (par rapport au côté du carré)** : Pour justifier ce que dit Aristote dans *Seconds Analytiques* concernant la démonstration de l'incommensurabilité de la diagonale dont le contraire impliquerait que 'tous les impairs deviendraient pairs', une preuve simple est donnée en comptant le nombre maximum de fois qu'un entier peut être divisé par *2*, c'est une généralisation du cas des pairs, au moins 1 fois divisible par *2*, et impairs qui ne sont pas divisibles par *2*. Je ne vais pas entrer dans le détail ce qui serait trop long, je renvoie à mon article de 2010 donné en bibliographe. Ce qui nous importe ici est qu'un corollaire quasi-immédiat de la démonstration montre que la démonstration d'incommensurabilité du côté d'un carré de surface entière par rapport à l'unité se ramène toujours au seul cas des surfaces impaires.

Plus précisément,

- Si un entier est divisible par *2* au plus un nombre impair de fois, alors le côté du carré de surface ce nombre est incommensurable à l'unité.
- Si au contraire, il est divisible par *2* un nombre pair de fois, alors il est commensurable à l'unité ou pas selon que le nombre obtenu après toutes les divisions par *2* est commensurable ou ne l'est pas. Et comme ce nombre est impair, on est ramené au cas des impairs.

    **Exemples** : *2* et *6* sont divisibles 1 fois par *2*, le côté des carrés de *2* et *6* pieds sont donc incommensurables à l'unité ; *8* est divisible 3 fois par *2* et donc là encore le côté du carré de surface *8* pieds est incommensurable à l'unité.
    *12* est divisible 2 fois par *2* (*12 = 4 × 3*), donc le côté du carré de *12* pieds est commensurable à l'unité si et seulement si celui du carré de surface *3* pied l'est. En termes modernes, $\sqrt{12}$ est rationnel si $\sqrt{3}$ l'est. Cela n'a rien de surprenant puisque $\sqrt{12} = \sqrt{(4 \times 3)} = 2\sqrt{3}$ !

En appliquant ce résultat à tous les entiers compris entre *2* et *17* (inclus ou pas), on peut éliminer tous les pairs. La suite étudiée par Théodore est donc formée des entiers impairs compris entre *3* et *17* (inclus ou pas) :

**Suite de Théodore** : *3, 5, 7, 9, 11, 13, 15, (17)*

Elle est plus courte (7 ou 8 termes) que la suite standard (contre 11 ou 12 termes) et va jouer un rôle important pour les preuves d'incommensurabilité que l'on considérera dans le second exposé.

### IV. La méthode de Théodore
#### i) Quelques méthodes antérieures

La méthode utilisée par Théodore pour construire la suite de ces dessins a été également beaucoup discutée. Certains pensent qu'elle n'importe pas puisque Platon ne la présente pas, d'autres considèrent la méthode dite par '*anthyphérèse*' qu'on peut exclure si l'on suit les contraintes données par Knorr, vu sa difficulté et sa longueur, pour une critique détaillée je vous renvoie à *EEE*, p. 118-126 (voir aussi son analyse sur la question des processus infinis pour les mathématiciens grecs, p. 36). Une autre est la construction directe par spirale assez populaire encore actuellement pour sa simplicité mathématique.

Je cite Szabó (*The Beginnings of Greek Mathematics*, p. 55-56) :

> 'Les commentateurs précédents ont invariablement demandé comment Théodore avait construit ces carrés ; (…) certains ont dit que cela n'avait pas d'importance, (…) ou alors ont adopté la tentative intéressante de reconstruction de la méthode donnée par Anderhub' :

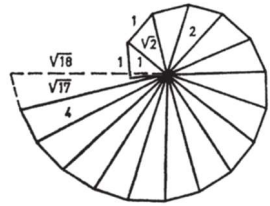

**Figure 2**

Szabó donne alors plusieurs raisons pour lesquelles cette construction n'est pas crédible. Nous pensons de même, mais d'abord pour une question textuelle : la construction ci-dessus nécessite l'étude de *tous* les carrés de surface compris entre 1 et 17 pieds, ce qui contredit directement Théétète, selon lequel Théodore commence avec celui de 3 pieds et saute celui de 4 pieds. Et il y a bien d'autres difficultés techniques pour la mise en œuvre concrète de cette méthode.

### ii) Notre reconstruction

La méthode qui suit est basée sur un corollaire immédiat du théorème de Pythagore obtenu simplement en itérant trois fois ce théorème (cf. Annexe I) :

**Corollaire du théorème de Pythagore**

Dans un triangle rectangle ABC, de hauteur BH, alors le carré de côté BH est égal au rectangle de côtés AH et HC (en écriture moderne $BH^2 = AH \times HC$) :

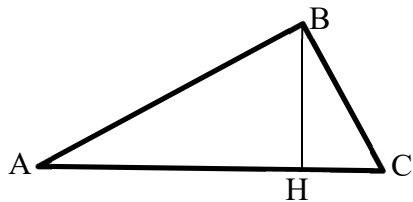

**Figure 3**

Elle utilise également un très ancien résultat attribué à Thalès : l'angle interceptant un diamètre est un angle droit (*Éléments* d'Euclide, proposition III.31). Dans la figure ci-dessous, tous les angles ADC, AEC et AFC sont droits.

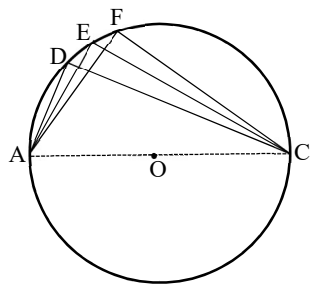

**Figure 4**

## V. Les dessins de Théodore

Tout d'abord Theodore trace une droite notée PZ en marquant des repères chaque pied, donc les points P, Q, A, B, C, D, etc. où PQ = QA = AB = … = 1 pied sur une longueur de 9 pieds. Puis il trace la droite perpendiculaire à PZ passant par Q, soit QT.

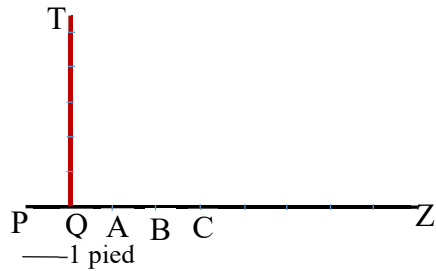

**Figure 5**

Puis viennent les différents dessins.

**Premier dessin (cas de 3 pieds)**

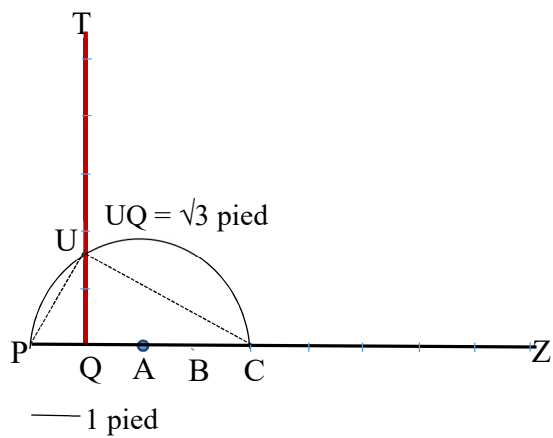

**Figure 6**

On considère d'abord le demi-cercle de centre A passant par P, donc de diamètre PC = 4 pieds. Puisque l'angle PUC est droit car intersectant un diamètre, d'après le corollaire du théorème de Pythagore, on a : $UQ^2 = PQ \times QC$. Mais PQ = 1 pied et QC = 3 pieds, donc en notations modernes $UQ = \sqrt{3}$ pieds (environ 0,53m).

**Deuxième dessin (cas de 5 pieds)**

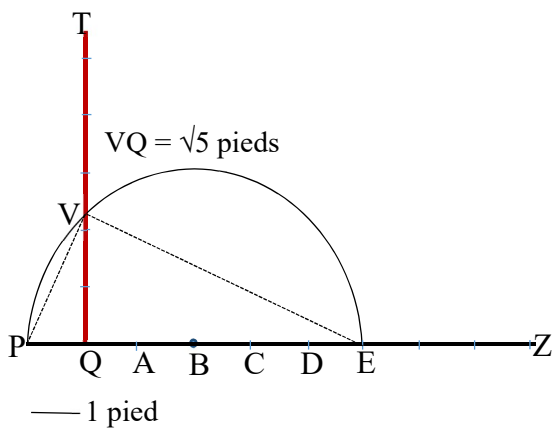

**Figure 7**

On décale alors le centre d'un pied en B et on considère le demi-cercle de centre B passant par P, donc de diamètre PC = 6 pieds. Cette fois c'est l'angle PVE qui est droit et puisque QE = 5 pieds, on obtient : $VQ^2 = PQ \times QE = 5$ pieds, d'où : $VQ = \sqrt{5}$ pieds (environ 0.65 m).

Théodore continue en décalant à chaque fois le centre du demi-cercle passant par P d'un pied. On obtient donc les demi-cercles de centres successifs, C, D, E, F, G, H, de diamètre respectifs 8, 10, 12, 14, 16 et 18 pieds. L'intersection de ces demi-cercles avec la perpendiculaire QT donne successivement les côtés de longueurs $\sqrt{7}$, $\sqrt{9}$, $\sqrt{11}$, $\sqrt{13}$, $\sqrt{15}$ et $\sqrt{17}$ pieds, inclus ou pas. Regardons ce qui se passe avec la dernière construction (éventuelle), celle de 17 pieds.

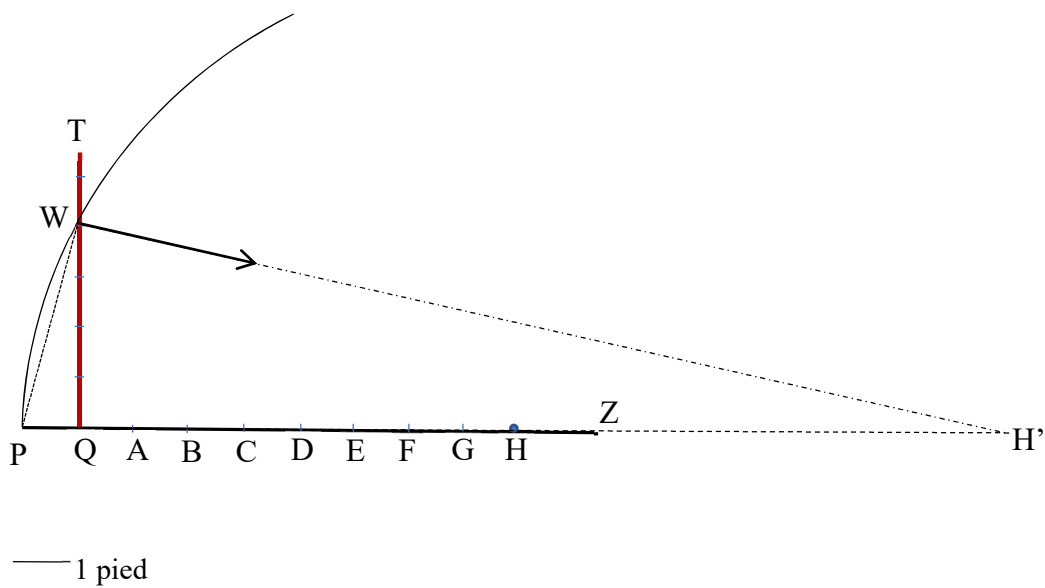

**Figure 8**

Le rayon PH = 9 pieds, le diamètre PH' = 18 pieds d'où QH' = 17 pieds. Soit W l'intersection du demi-cercle avec la perpendiculaire QT. L'angle PWH' intersectant un diamètre est droit, et d'après le corollaire, on a donc : $QW^2 = PQ \times QH'$, d'où QW = $\sqrt{17}$ pieds (environ 1.25 mètre).

On a un témoignage historique allant dans le sens de cette méthode que l'on retrouve dans ce qu'on appelle le *Commentaire Anonyme* (cf. Bastiani-Sedley donné en bibliographie).

### iii) Le dessin final

On obtient alors comme dessin final :

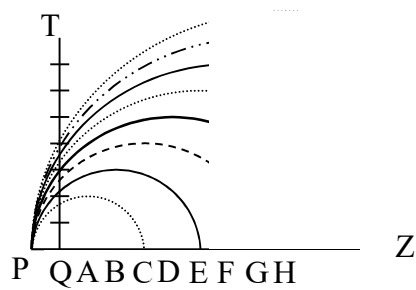

**Figure 9**

**Conclusion.** Cette construction est simple et facile à exécuter, contrairement à celle de l'*anthyphérèse*. Elle permet de tracer les diverses grandeurs incommensurables étudiées par Théodore au cas par cas, à tour de rôle ('κατὰ μίαν ἑκάστην'), sans bouger sans cesse ni tracer des perpendiculaires dans toutes sortes de directions, comme dans le cas de la spirale. Elle est donc un bon candidat pour la construction des longueurs incommensurables de Théodore. En fait, elle est compatible avec les 9 conditions imposées, et c'est la seule. Je ne vais pas le montrer ici en détails, car ces contraintes sont surtout importantes pour les preuves utilisées par Théodore pour ses démonstrations d'incommensurabilité, et d'ailleurs c'est par rapport à ces elles que Knorr a formulé les siennes. Ce sera donc dans le second exposé qu'elles seront considérées, et qu'on montrera comment et pourquoi notre reconstruction est bien la seule qui les vérifie.

J'ai considéré ici la partie graphique de la leçon de Théodore, car Théodore est présenté d'abord comme un géomètre, bien que le terme lui-même était à l'époque très général et souvent synonyme de mathématicien. Inversement, les mathématiques de Platon sont d'abord fondées sur des questions de théorie des nombres, même lorsqu'il donne des exemples géométriques comme ici ou dans le *Ménon*. On peut même aller jusqu'à dire que pour Platon, la géométrie serait une sorte de branche 'd'une arithmétique supérieure' (*EEE*, p. 90)).

Nous en verrons un exemple dans le deuxième exposé la semaine prochaine.

**Annexe I**

**Corollaire du corollaire théorème de Pythagore** : *Dans un triangle rectangle ABC, de hauteur BH, alors le carré de côté BH est égal au rectangle de côtés AH et HC (en écriture moderne BH² = AH×HC).*

**Démonstration.** Pour la clarté et la brièveté, on utilise ici les notations algébriques modernes, le passage dans le langage de la géométrie grecque ancienne étant aisé.

Soit ABC un triangle rectangle en B, et BH la hauteur issue de B :

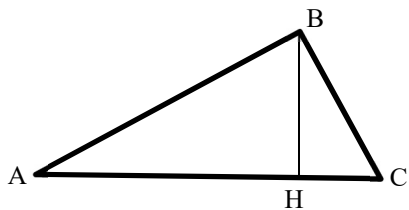

**Figure 10**

Il faut montrer que : AB× BC= BH²

(ou en termes plus proches de la géométriques grecque ancienne : le rectangle de côtés AB, BC est égal au carré de côté BG=BH).

La démonstration se fait par applications successives du théorème de Pythagore aux triangles rectangles BHA, BHC et ABC :

AB² = AH² + HB² (Pythagore pour BHA)

BC² = HC² + BH² (Pythagore pour BHC)

AC² = AB² + BC² (Pythagore pour ABC).

Puisque : AC² = (AH + HC)²= AH² + HC² + 2AH×HC (c'est la proposition II.4 des *Éléments*) , la dernière égalité donne :

AH² + HC² + 2AH×HC = AB² + BC².

On a donc d'après les deux premières égalités :

AH² + HB² + 2AH×HC = AH² + HB² + HC² + BH², d'où :

2AH×HC = 2 HB², donc : AH×HC = 2 HB². CQFD

Ce résultat est le corollaire de la proposition VI.8 des *Éléments* d'Euclide, mais la démonstration donnée ici est différente.

**Annexe II : Les propriétés nécessaires pour une reconstitution acceptable des démonstrations d'incommensurabilité d'après Wilbur Knorr (EEE, p. 96-97)**

'From our above examination of the Theaetetus-passage we may deduce the following set of requirements which must be satisfied by any acceptable mathematical reconstruction of Theodorus' proofs:

(a) The proofs are demonstratively valid.
(b) The treatment by special cases and the stopping at 17 are necessitated by the methods of proof employed.
(c) The proofs will be understood to apply to an infinite number of cases.
(d) No use may be made of the dichotomy or square and oblong numbers in Theodorus' studies, either in the demonstrations or in the choice of cases to be treated.
(e) Theodorus' proofs utilize the special relations of the lines drawn in the construction of the *dynameis*. The geometrical methods of construction are of the type characteristic of metrical geometry as developed in *Elements* II and are closely associated with a certain early style of arithmetic theory.
(f) But the arithmetic methods by which Theaetetus could prove the two general theorems, on the incommensurability of lines associated with non-square and non-cubic integers, were not available to Theodorus.

These follow from a strict reading of the passage. If a given reconstruction satisfies them, we may accept it as possible; otherwise, not. It is conceivable, of course, that any number of reconstructions will be possible in answer to a finite set of conditions. But we will find that, in fact, no reconstruction (prior to the present one) satisfies them all.'

# La leçon mathématique de Théodore (*Théétète*, 147c7-d6) II

## Les démonstrations d'incommensurabilité de Théodore

Salomon Ofman

## VI. Introduction

Ce deuxième exposé porte sur les méthodes utilisées par Théodore dans ses démonstrations d'incommensurabilité. Il est basé comme premier sur un travail en commun avec Luc Brisson sur le passage mathématique et son articulation au reste du *Théétète*.

Dans cet exposé, j'utiliserai les notations modernes pour simplifier les démonstrations mathématiques. Dans les annexes, vous trouverez les démonstrations écrites dans le langage des mathématiques grecques anciennes.

## VII. Le texte

| **Thééthète** | **Θεαίτητος** |
|---|---|
| Et en fait, il se pourrait que tu demandes la sorte de chose qui nous est venue à l'esprit, à moi et à ce Socrate-là ton homonyme, alors que nous discutions entre nous tout à l'heure. | (147c7) ἀτὰρ κινδυνεύεις ἐρωτᾶν οἷον καὶ αὐτοῖς ἡμῖν ἔναγχος εἰσῆλθε διαλεγομένοις, ἐμοί τε καὶ τῷ σῷ ὁμωνύμῳ τούτῳ Σωκράτει. |
| **Socrate** | **Σωκράτης** |
| Quoi donc, Théétète ? | τὸ ποῖον δή, ὦ Θεαίτητε; |
| **Thééthète** | **Θεαίτητος** |
| Quelque chose à propos des puissances. Théodore ici présent nous avait trace des figures, mettant en évidence que la puissance de 3 pieds aussi bien que celle de 5 pieds ne sont pas commensurables en longueur avec celle de 1 pied ; et il a continué ainsi, prenant chacune à tour de rôle, jusqu'à celle de 17 pieds. Il s'est arrêté là pour quelque raison. | (147d2) περὶ δυνάμεων τι. ἡμῖν Θεόδωρος ὅδε ἔγραφε, τῆς τε τρίποδος πέρι καὶ πεντέποδος ἀποφαίνων ὅτι μήκει οὐ σύμμετροι τῇ ποδιαίᾳ, καὶ οὕτω κατὰ μίαν ἑκάστην προαιρούμενος μέχρι τῆς ἑπτακαιδεκάποδος· ἐν δὲ ταύτῃ πως ἐνέσχετο. |

Comme je vous l'ai dit la semaine dernière, ce passage a donné lieu a beaucoup de reconstructions divergentes, la meilleure étant sans doute celle de Wilbur Knorr dans son ouvrage *The Evolution of Euclidean* Elements (ou *EEE*), bien que de nombreux points posent de grandes difficultés. Ce qui nous paraît important, c'est que contrairement aux autres commentateurs, il ne se contente pas de détailler la méthode qu'il juge la meilleure. Il donne six conditions que toute reconstruction de la méthode de Théodore doit vérifier pour être acceptable, ouvrant la voie théoriquement à une indéfinité de reconstructions possibles différentes et opposées à la sienne. A ces 6 conditions, nous en avons ajouté trois autres. Nous listons toutes ces neuf conditions ci-dessous.

## VIII. Les six conditions de Knorr

Nous traduisons le texte de *EEE*.

'De notre examen ci-dessus du passage du *Théétète*, nous pouvons déduire l'ensemble des exigences suivantes que doit vérifier toute reconstruction mathématique des preuves de Théodore pour être acceptable :

1) Les preuves doivent être correctes.

2) Le traitement par cas particuliers et l'arrêt en 17 sont nécessités par les méthodes de preuve employée.

3) Les preuves doivent être comprises comme pouvant s'appliquer à un nombre infini de cas.

4) Aucune utilisation ne doit être faite de la dichotomie entre entiers carrés parfaits et entiers rectangulaires non carrés dans les études de Théodore, que ce soit dans les démonstrations ou dans les cas traités.

5) Les preuves de Théodore utilisent les relations spéciales entre les lignes tracées dans la construction des *dynameis*. Les méthodes géométriques de construction sont de type caractéristique de la géométrie métrique telle qu'elle est développée dans les *Éléments* II et sont étroitement associées avec un certain style du début de la théorie arithmétique.

6) Mais les méthodes arithmétiques par lesquelles Théétète pouvait prouver les deux théorèmes généraux concernant l'incommensurabilité des lignes associées à des entiers non-carrés et non-cubes parfaits n'étaient pas disponibles à Théodore.

Ceci suit d'une lecture stricte du passage. Si une reconstruction donnée les satisfait, on peut l'accepter comme possible, sinon non. Il est concevable naturellement qu'un nombre quelconque de reconstructions soient possibles en réponse à un ensemble fini de conditions. Mais nous nous apercevrons qu'en fait, aucune reconstruction (avant la présente [du *EEE*]) ne les satisfait toutes.' (*EEE*, p. 96-97).

## IX. Nos trois conditions supplémentaires

Comme je l'ai dit dans le premier exposé, non contents d'accepter le défi de ces six conditions, nous en ajoutons trois autres qui nous paraissent essentielles, et que la reconstruction de Knorr, et tout aussi bien les autres méthodes, ne vérifie pas :

4) On ne doit pas postuler *a priori* que Platon veut célébrer ici Théétète, futur grand mathématicien, et qu'il lui rend hommage. Il faut juger sur les textes.

5) Le passage doit être lu comme partie du *Théétète* et non pas comme une partie autonome.

6) La description par Platon est réaliste au sens où elle doit suivre de près les leçons réelles des mathématiciens sur le sujet, que cette leçon soit mise en scène par Platon ou ait vraiment eu lieu.

Dans l'exposé précédent, j'ai montré les graphiques tracés par Théodore, et la manière dont il les construisait. Dans celui-ci, je vais considérer sa méthode pour obtenir l'incommensurabilité de certaines grandeurs par rapport au pied-unité. Cette reconstruction vérifie ces neuf conditions.

## X. Les reconstructions

De manière générale, toutes les reconstructions commencent à la manière de la démonstration d'incommensurabilité de la diagonale par rapport au côté du carré. La méthode est une démonstration par l'absurde ou l'impossible : on suppose donc que le côté du carré est commensurable pour aboutir à une conclusion absurde ou impossible.

**Deux reconstructions usuelles**

Deux grands types de méthodes sont communément admises parmi les commentateurs, avec toutefois de nombreuses variantes suivant les auteurs.

1) La première méthode qui a eu longtemps les faveurs des spécialistes est très proche de celle-là même qu'on utiliserait encore aujourd'hui. Elle s'appuie sur une proposition des *Éléments* d'Euclide démontrant l'incommensurabilité de la diagonale du carré, une proposition notée X. 117. Cette proposition est considérée à juste titre comme inauthentique car elle n'est pas cohérente avec le reste du traité d'Euclide traitant de ces questions. Je vous renvoie pour les détails à mon article de 2010 donné en bibliographie.

   Mais qu'en est-il pour les démonstrations de Théodore ? Même inauthentique, ce pourrait être sa méthode de preuve. Cela se heurte à deux problèmes essentiels. Cette méthode suppose un développement de la théorie des nombres qui n'était certainement pas disponible à l'époque de Théodore. Et surtout, elle est en contradiction avec le texte car elle donne une solution générale au problème d'incommensurabilité, alors que Théétète précise que Théodore procède au cas par cas, à tour de rôle ('κατὰ μίαν ἑκάστην'). En outre, elle n'explique pas l'arrêt à *17*. En résumé, elle viole la condition 2 de Knorr.

2) La seconde méthode est celle dite par '*anthyphérèse*' ou 'soustraction alternée' déjà évoquée à l'exposé précédent. Elle a eu une grande popularité parmi les chercheurs du milieu du 20ème siècle, au point de concurrencer la précédente, et elle a encore des partisans aujourd'hui. Brièvement, c'est ce qu'on appelle en langage moderne la théorie des fractions continues. Elle s'appuie sur le début du livre X des *Éléments* d'Euclide où l'*anthyphérèse* de 2 grandeurs est définie (proposition 1) et où il est montré que si le processus de 'soustractions alternées' ne s'arrête pas, alors les 2 grandeurs sont incommensurables (proposition 2).
Théodore l'utiliserait pour démontrer les cas d'incommensurabilité des côtés des carrés. J'ai déjà dit dans l'exposé précédent que cela est peu vraisemblable en raison de la difficulté et de la longueur des démonstrations, même si on a proposé toutes sortes de subtilités mathématiques pour la simplifier. Elle viole par ailleurs plusieurs des 9 conditions que nous avons posées (particulièrement les conditions 2, 3, 5 et surtout 9). Pour une critique détaillée portant sur les conditions 2, 3 et 5, je vous renvoie à Knorr, *EEE* p. 118-126 ainsi qu'à son analyse sur la question des processus infinis pour les mathématiciens grecs, p. 36.

On peut alors se demander pourquoi elle a été acceptée si facilement. C'est que jusqu'à une cinquantaine d'années, il n'y avait d'autre choix que la méthode basée sur la proposition inauthentique X.117. Et pour reprendre les mots de Socrate, cette dernière est 'encore plus impossible' (192b6) que celle par '*anthyphérèse*' dont les démonstrations d'incommensurabilité sont faites effectivement au moins au cas par cas, comme l'exige le texte du *Théétète*.

XI. **Notre reconstruction**

1) **Préliminaires**

   Elle s'appuie sur un résultat sur les carrés impairs dont le premier à avoir souligné l'intérêt pour les preuves d'incommensurabilité a été Jean Itard dans son ouvrage sur les livres arithmétiques d'Euclide que vous trouverez en bibliographie. Il s'agit d'un résultat très simple concernant les carrés impairs, qui apparaît facilement dès que l'on considère des tableaux de carrés, de tels tableaux remontant aux Babyloniens. C'est sans doute la raison pour laquelle Itard se contente d'en donner une preuve sous une forme algébrique moderne. Knorr la redémontre sous la forme d'un 'théorème 10' et en donne une preuve dans le style dit-il de l'arithmétique grecque ancienne. Il apparaît sous une forme ou une autre dans plusieurs textes de l'Antiquité qui nous sont parvenus, tels l'*Expositio* de Théon, l'*Arithmetica* de Diophante, *In Nichomachum* de Jamblique, mais également dans les *Platonicae Questiones* chez Plutarque. Leur analyse ainsi que les méthodes utilisées permettent à Knorr de conclure que c'était un théorème connu des Pythagoriciens anciens (*EEE*, p. 154).

## 2) Un résultat sur les carrés impairs

Soit alors un nombre impair. C'est un multiple de *2* augmenté d'une unité. Si on regarde la suite des premiers carrés impairs, on a :

*$3^2$ = 9 = 8+1 ; $5^2$ = 25 = 24+1 ; $7^2$ = 49 = 48+1 ; $9^2$ = 81 = 80+1 ; $11^2$ = 121 = 120+1 ; $13^2$ = 169 = 168+1 …,*

et tous ces nombres pairs qui apparaissent alors : *8, 24, 48, 80, 120, 168* sont non seulement pairs, non seulement des multiples de *4*, mais des multiples de *8*. Cela est vrai en général, et on a (cf. Annexe II) :

**Résultat des carrés impairs** : *Le carré d'un impair est un multiple de 8 augmenté d'une unité.*

## 3) Les démonstrations d'incommensurabilité

Le raisonnement utilise les multiples de 8, c'est-à-dire :

*8 ; 16 ; 24 ; …*

Pour les modernes, le premier multiple serait 0 qui n'est pas considéré comme un nombre dans les mathématiques grecques anciennes et n'apparaît donc pas ici.

### 1. Le cas du carré de *3* pieds

Nous avons vu dans l'exposé précédent comment Théodore dessinait le côté du carré. Maintenant il va montrer que ce côté est incommensurable à l'unité. Cela se fait par un raisonnement sur les entiers.

En notations modernes, on suppose donc que $\sqrt{3}$ pieds est commensurable à l'unité, donc il existe des entiers *p* et *q* tels que

$\sqrt{3}/1 = p/q$

d'où :

$\sqrt{3} \times q = p \times 1 = p$

et en prenant leurs carrés :

*$3q^2 = p^2$.*

**Remarque.** On ne suppose pas que *p* et *q* sont premiers entre eux comme cela est fait dans les méthodes usuelles, par contre le comptage du nombre de divisions possibles par *2* à la base de la démonstration de l'incommensurabilité de la diagonale par rapport au côté du carré, permet de prendre de *p* et *q* tous deux impairs.

D'après le résultat des carrés impairs ci-dessus, $q^2$ est un multiple de *8* augmenté d'une unité. Donc $3q^2$ est un multiple de *8* augmenté de 3 ($3q^2 = 8k+3$). Mais $p^2$ est un multiple de *8* augmenté d'une unité ($p^2 = 8h+1$), donc puisqu'il est égal à $3q^2$ il est aussi un multiple de $8$ augmenté de $3$. On a donc :

$8k+3 = 8h+1$,

d'où :

$2 = 8h-8k$.

Mais la différence de 2 multiples de 8 est un multiple de 8, ce qui est impossible car *2* n'est pas un multiple de *8*.

Conclusion : *le côté du carré de (surface) 3 pieds est incommensurable à 1 pied.*

### 2. Le cas du carré de *5* pieds

On procède de même et on obtient :

$\sqrt{5}/1 = p/q$ d'où $\sqrt{5} \times q = p \times 1 = p$

d'où en prenant leurs carrés :

$5q^2 = p^2$ avec *p* et *q* impairs.

De même que ci-dessus, d'après le résultat sur les carrés impairs, $p^2$ et $q^2$ sont des multiples de *8* plus une unité, et en particulier : $5q^2 = 8k+5$. On a donc :

$5q^2 = 8k+5 = p^2 = 8h+1$, d'où $8k+5 = 8h+1$, et donc :

$8k+4 = 8h$, d'où $4 = 8h-8k$.

Mais alors *4* serait la différence de 2 multiples de 8, ce qui est impossible car *4* n'est pas un multiple de *8*.

Conclusion : *le côté du carré de (surface) 5 pieds est incommensurable à 1 pied.*

### 3. Le cas du carré de *7* pieds

De même, on obtient :

$\sqrt{7}/1 = p/q$ d'où $\sqrt{7} \times q = p \times 1 = p$

d'où en prenant leurs carrés :

$7q^2 = p^2$ avec *p* et *q* impairs.

D'où :

$7q^2 = 8k+7 = p^2 = 8h+1$, d'où $8k+7 = 8h+1$, et donc $6 = 8h-8k$,

ce qui est impossible puisque *6* n'est pas un multiple de *8*.

Conclusion : *le côté du carré de (surface) 7 pieds est incommensurable à 1 pied.*

### 4. Le cas du carré de *9* pieds

Dans ce cas, on peut écrire :

$\sqrt{9}/1 = p/q$ d'où $\sqrt{9} \times q = p \times 1 = p$

d'où en prenant leurs carrés :

$9q^2 = p^2$ avec *p* et *q* impairs.

D'où :

$9q^2 = 8k+9 = p^2 = 8h+1$

d'où :

$8k+9 = 8h+1$, et donc $8 = 8h-8k$.

Mais cela est possible puisque $8$ est un multiple de $8$ !

Or $9$ est un carré parfait, et le côté de longueur $\sqrt{9}$ pieds = $3$ pieds qui est commensurable à $1$ pied.

## 5. Le cas du carré de *11* pieds

De même, on obtient :

$\sqrt{11}/1 = p/q$ d'où $\sqrt{11} \times q = p \times 1 = p$

d'où en prenant leurs carrés :

$11q^2 = p^2$ avec $p$ et $q$ impairs.

D'où :

$11q^2 = 8k+11 = p^2 = 8h+1$,

d'où :

$8k+11 = 8h+1$,

et donc :

$10 = 8h-8k$.

Ce qui est impossible puisque $10$ n'est pas un multiple de $8$.

Conclusion : le côté du carré de (surface) *11* pieds est incommensurable à *1* pied.

## 6. Le cas du carré de *13* pieds

On a :

$\sqrt{13}/1 = p/q$ d'où $\sqrt{13} \times q = p \times 1 = p$, et donc : $13q^2 = p^2$ avec $p$ et $q$ impairs.

D'où :

$13q^2 = 8k+13 = p^2 = 8h+1$, d'où $8k+12 = 8h$, et donc $12 = 8h-8k$.

Ce qui est impossible puisque $12$ n'est pas un multiple de $8$.

Conclusion : le côté du carré de (surface) *13* pieds est incommensurable à *1* pied.

## 7. Le cas du carré de *15* pieds

On a :

$\sqrt{15}/1 = p/q$ d'où $\sqrt{15} \times q = p \times 1 = p$, et donc : $15q^2 = p^2$ avec $p$ et $q$ impairs.

D'où :

$15q^2 = 8k+15 = p^2 = 8h+1$, d'où $8k+14 = 8h$, et donc $14 = 8h-8k$.

Ce qui est impossible puisque $14$ n'est pas un multiple de $8$.

Conclusion : le côté du carré de (surface) *15* pieds est incommensurable à *1* pied.

## 8. Le cas du carré de *17* pieds

On a encore :

$\sqrt{17}/1 = p/q$ d'où $\sqrt{17} \times q = p \times 1 = p$, et donc : $17q^2 = p^2$ avec *p* et *q* impairs.

D'où :

$17q^2 = 8k+17 = p^2 = 8h+1$, d'où $8k+16 = 8h$, et donc $16 = 8h-8k$.

Mais *16* est un multiple de *8*, et donc le côté du carré (de surface) *17* pieds pourrait être commensurable à l'unité. Mais la propriété des impairs est nécessaire mais non pas suffisante, on ne peut donc rien dire du moins par la méthode des carrés impairs, et Théodore est donc obligé de s'arrêter là.

La méthode de Théodore ne lui permet pas de dire si le côté d'un carré de surface *8k+1* qui n'est pas un carré parfait est ou n'est pas commensurable à *1* pied. Par exemple, il ne peut résoudre la question pour les carrés de surface *8×2 + 1 = 17* pieds, *8×4 + 1 = 33* pieds, *8×5 + 1 = 41* pieds, *8×7 + 1 = 57* pieds, *8×8 + 1 = 65* pieds, … Les seuls cas auxquels il peut répondre sont pour les carrés de surface *9 (= 3×3)* pieds, *25 (= 5×5)* pieds, *49 (= 7×7)* pieds, … c'est-à-dire pour les carrés parfaits.

**Vérification.** On va voir maintenant que cette méthode vérifie les 9 (6+3) conditions que nous avons données, et qu'il faut imposer à toute reconstruction pour être valide.

1) Les preuves doivent être correctes.
   *C'est bien le cas.*
2) Le traitement par cas particuliers et l'arrêt en 17 sont nécessités par les méthodes de preuve employée.
   *C'est bien le cas.*
3) Les preuves doivent être comprise comme pouvant s'appliquer à un nombre infini de cas.
   *C'est bien le cas.*
4) Aucune utilisation ne doit être faite de la dichotomie entre entiers carrés parfaits et entiers rectangularies non carrés dans les études de Théodore, que ce soit dans les démonstrations ou dans les cas traités.
   *C'est bien le cas. En fait, contrairement à Knorr lui-même, le carré de surface 9 un carré parfait est dans la liste des grandeurs considérés par Théodore selon notre reconstruction.*
5) Les preuves de Théodore utilisent les relations spéciales entre les lignes tracées dans la construction des *dynameis*. Les méthodes géométriques de construction sont de type caractéristique de la géométrie métrique telle qu'elle est développée dans les *Éléments* II et sont étroitement associées avec un certain style du début de la théorie arithmétique.
   *Là encore c'est le cas, la géométrie métrique est celle qui travaille sur les mesures et les nombres, ce qui est le cas ici. D'ailleurs Knorr utilise essentiellement la propriété des carrés impairs dans sa reconstruction.*
6) Mais les méthodes arithmétiques par lesquelles Théétète pouvait prouver les deux théorèmes généraux concernant l'incommensurabilité des lignes associées à des entiers non-carrés et non-cubes parfaits n'étaient pas disponibles à Théodore.
   *C'est bien évident, puisque la méthode ne permet de rien dire concernant une infinité de carrés de surfaces nombres entiers, à savoir tous ceux de la forme 8k+1 qui ne sont pas des carrés parfaits.*
7) On ne doit pas postuler *a priori* que Platon veut célébrer Théétète, futur grand mathématicien, et qu'il lui rend hommage ici.

*Absolument. Au vu de la solution donnée par la méthode de Théodore, on peut d'ailleurs mettre également en question l'admiration que Platon portait pour celui-ci qui aurait été son « maître » en géométrie.*

8) Le passage doit être lu comme partie du *Théétète* et non pas comme une partie autonome. *Absolument. C'est pourquoi nous commençons avant ce qu'ont fait les auteurs qui ont travaillé sur ce passage pour l'introduire dans son contexte. Il y a bien d'autres points importants dont nous n'avons pas parlés et qui relient ce passage à la question de l'erreur, de l'opinion, de l'opinion vraie face à la science, mais aussi les différentes définitions de la science, dont les trois dernières concernant le* logos.

9) La description par Platon est réaliste au sens où elle doit suivre de près les leçons réelles des mathématiciens sur le sujet, que cette leçon soit mise en scène par Platon ou ait eu vraiment lieu. *C'est pourquoi j'ai donné la totalité des cas d'incommensurabilités de la leçon de Théodore. Le point faible des autres reconstructions est qu'elles doivent considérer que Platon symbolise ici, car dans toutes les autres méthodes proposées, il est totalement irréaliste que l'on puisse montrer ces incommensurabilités dans le cadre d'une seule leçon de mathématique. Cela est valable pour la méthode longue et compliquée donnée par Knorr comme il le reconnaît d'ailleurs, ajoutant toutefois que cela n'importe pas vraiment car, dit-il, les reconstructions antérieures à la sienne sont encore pires, en quoi il a raison.*

**Conclusion.** La méthode présentée ici utilise des résultats élémentaires et anciens fondés essentiellement sur les questions du pair et de l'impair. Théodore en suivant cette méthode a pu la donner *verbatim* sans qu'on ait à supposer que Platon fait dans le symbolisme ou la métaphore. C'est ce à quoi doivent se résoudre les autres reconstructions proposées par les nombreux commentateurs et historiens des mathématiques, y compris en fin de compte celle de Knorr. Cela est vrai, que cette leçon de Théodore à Théétète et ses jeunes camarades, ait eu lieu ou soit entièrement l'œuvre de l'imagination de Platon. C'est également la seule qui soit entièrement cohérente avec le texte platonicien. Pour terminer, je voudrais revenir sur ce qui est sans doute l'aspect principal pour Platon du passage mathématique, et en particulier de la leçon de Théodore, c'est-à-dire son intérêt philosophique.

Pour cela, je vais très brièvement considérer la dernière définition de la science par Théétète, à savoir la *doxa* vraie à laquelle s'ajoute un *logos* (201c9-d1), puis les trois définitions du *logos* (que plusieurs interprètes de Platon remplacent par 'logismos', 'raisonnement' suivant ce que dit Socrate dans le *Ménon*), ce qui lui permet de réfuter une dernière fois Théétète (202d-210b). En dépit de cette réfutation, certains commentateurs tels Christopher Rowe qui par ailleurs a donné une très bonne traduction du *Théétète*, voient dans cette définition, celle de la science pour Platon. Nous pensons que la leçon de mathématique telle que nous l'avons considérée, entre autres choses, va directement à l'encontre de cette thèse.

Tout d'abord, cette leçon est incontestablement un discours vrai, le raisonnement est correct et le résultat est juste. Pourtant, il est peu probable que Platon s'en contenterait et le considérerait comme scientifique.

En effet, d'une part il ne donne qu'un résultat partiel, sans répondre à la question générale 'quand le côté d'un carré entier est-il commensurable à l'unité ?'. Ce résultat est obtenu via un résultat sur les carrés impairs qui n'a que peu, voire pas du tout de relation avec la question posée. En outre, la réponse étant partielle, même si on obtient ainsi une solution pour plus de 80% des cas, elle ne peut prétendre donner la cause de la commensurabilité ou de l'incommensurabilité des côtés des carrés entiers. Or dans le *Ménon* encore, Socrate exige d'un discours scientifique qu'il donne la cause de ce sur quoi il porte, qu'il le lie par un 'raisonnement sur la cause' ('αἰτίας λογισμῷ') (*Ménon*, 98a3-4). Il n'est pas même nécessaire de poser le problème délicat des deux types de discours dont l'un porterait sur les réalités intelligibles et l'autre sur les choses sensibles, établissant ainsi une séparation radicale entre *doxa* et science. La dernière définition du *Théétète* est déjà réfutée par la leçon mathématique, et cela quelle que soit la définition que l'on donne au *logos*.

Dans cette très brève analyse de conclusion j'ai voulu présenter de manière très partielle l'aspect central que tient selon nous le passage mathématique, non pas seulement du point de vue de l'histoire des mathématiques pré-euclidiennes, mais pour l'interprétation philosophique du *Théétète* et de son enquête sur une définition de la science, mais aussi de la *doxa*, du raisonnement (*logismos*) et du discours (*logos*).

# Annexe I

'the following general features of the mathematics Plato introduces. First, the examples cited are most usually drawn from the theory of number, and many of the geometric examples are in fact very closely associated with arithmetic studies. Indeed, Plato sometimes even speaks as if the fields of geometry and stereometry were branches of higher arithmetic. But this would be an exaggeration. What we see, rather, is that Plato's instances of geometry always fall into the pattern of metrical geometry, and in particular, of a metrical geometry directed tow. and the elucidation of the properties of integers. The overwhelming majority of mathematical examples in Plato's works are arithmetical. His favorite example is that of the odd and even. Among other arithmetical examples are the 'nuptial number' and the 'rational and irrational diameters' *(Republic* 546C) and the 'tyrant's number' *(Republic* 587E). Even when Plato introduces materials we would properly designate as geometric, the arithmetic intent is present. At *Meno* 82-85, for instance, the subject of discussion is a geometric theorem: that the double-square is constructed as the square on the diagonal of the given square. But the assignment of the 'foot' as unit establishes an affinity with the tradition of metrics and the mathematics of the passage is executed entirely within the field of arithmetic: only the operations of counting and of multiplication are employed in the analysis. When we consider the practice of representing numbers by geometric figures (the number two, for instance, may be figured as a rectangle of sides one and two units), we may view the Meno-passage as an investigation of the nature of such representations: how can the number eight be constructed as a square? We see the relevance of the concept of *dynamis* here, made explicit at *Statesman* 266A, and of course developed fully at *Theaetetus* 147D.' La terminologie de Théodore est métrique. Il utilise l'unité concrète du pied plutôt que les unités abstraites tel que le fait Euclide. (*EEE*, p. 90).

**Annexe II**

**Démonstration du résultat sur les carrés impairs**. *Tout carré impair est un multiple de 8 augmenté d'une unité.*

**Démonstration 1 (algébrique moderne).** On va tout d'abord donner une démonstration en utilisant le symbolisme moderne :

$m^2 = (2k+1)^2 = 4k^2 + 4k + 1 = 4(k^2 + k) + 1 = 4k(k+1) + 1$.

Comme ou bien $k$ ou bien $k+1$ est pair, $m^2$ est un multiple de *8* augmenté de *1*. En écriture en mathématiques grecques anciennes, on aurait : le carré construit (= de côté) *m* est égal à celui construit sur *2k+1*, donc (*Éléments*, ???).

Nous allons également donner plusieurs démonstrations plus en ligne avec la mathématique grecque ancienne.

**Démonstration 2 (arithmétique grecque ancienne).** Cela est fait via un résultat sur les entiers pairs.

Tout d'abord on a :

i) On montre tout d'abord que : *tout entier pair est ou bien lui-même ou bien si on lui ajoute 2 mesuré par 4.*
Soit *A* un tel entier et *B* ce nombre ajouté de *2*. Si *4* mesure *A*, le résultat est vrai. Soit donc *4* ne mesurant pas 4, c'est-à-dire il est impairement pair. Alors *A* est deux fois un impair *C*, donc la somme d'un pair *D* et de l'unité. Donc *A* est 2 fois *D* plus *2*, et *B* est 2 fois *D* plus *4*. Mais *D* est pair donc *4* mesure *2D*, donc il mesure *B*. [1]

ii) Soit alors *A* entier impair. C'est donc un pair plus une unité. Soit *B* cet entier pair. Le carré de *A* est égal à une unité augmentée du carré de *B* plus 2 fois *B*. Soit *C* ce nombre, il est égal au produit de *B* par *B* augmenté de *2*. D'après le résultat précédent, ou bien *4* mesure *B* et donc il mesure le carré de *B* et *C*, ou bien il mesure *2* ajouté à *B*. Soit *D* ce nombre. Mais *C* est le produit de *B* par *D* donc puisque *B* est pair, il est mesuré par *8*. [2] CQFD

**Démonstration 3** (Knorr, *EEE*, p. 153)

Il est donné sous la forme d'un théorème :

Théorème 10. *Tout nombre impair carré diminué d'une unité est 8 fois le multiple d'un nombre triangulaire.*

Knorr montre tout d'abord un autre théorème (*EEE*, p. :

Théorème 6 : *Tout nombre 'oblong' ['ἑτερομήκης', i.e. le produit de 2 entiers consécutifs] est le double d'un nombre triangulaire.*

**Démonstration du théorème 6**. Cela résulte immédiatement de la figure ci-dessous :

---

[1] Algébriquement, $B = A + 2 = 2C + 2 = 2(D+1) + 2 = 2D + 4$ avec *D* entier pair.
[2] Algébriquement, $A = B + 1$ and $(B+1)^2 = B(B+2) + 1$. Puisque *B* est pair, ou bien *B* ou *B+2* est un multiple de *4* (premier résultat), donc un multiple de *8*.

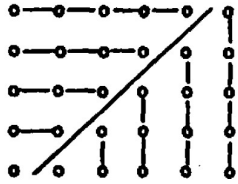

Cela peut être démontré en appliquant le Théorème 6 à la figure suivante :

Le théorème 10 résulte alors de l'application du théorème 6 à la figure suivante :

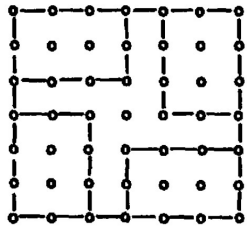

En effet, chaque portion 'oblongue' du carré impair dessiné dans cette figure peut être divisé en 2 nombres triangulaires égaux. Ainsi le grand carré, moins une unité, est égal à 8 nombres triangulaires égaux, comme l'affirme le théorème 10.

**Démonstration 4**

Analogue à la précédente mais plus brève car évident sur la figure :

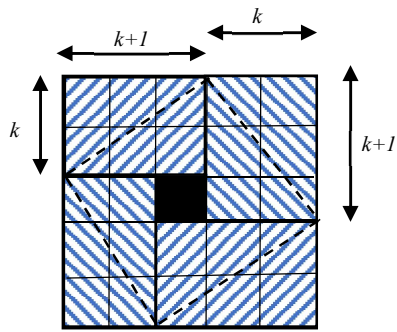

**Figure 11**

Le carré de côté *2k+1 = k + (k+1)* est somme d'une unité (le carré en noir) et de 4 rectangles de côtés (*k,k+1*). Le produit de 2 nombres consécutifs étant pair (appelé parfois '*hétéromèque*'), la surface du carré est égale à une unité plus un multiple de *8*. CQFD

# Bibliographie